\title{$S_3 \times S_3$ acts freely on $S^3 \times S^3$}
\author{James F. Davis\thanks{Partially support by an NSF grant.}\\
Department of Mathematics\\
Indiana University \\
Bloomington, IN 47405 \\
email: jfdavis@indiana.edu}
\date{}
\newcommand{\F}{\mathbb F}
\newcommand{\Q}{\mathbb Q}
\documentclass[12pt]{article}
\usepackage{amssymb,amsmath,amscd,amsthm}
\newcommand{\R}{\mathbb R}
\newcommand{\Z}{\mathbb Z}
\newtheorem*{theorem}{Theorem}
\newtheorem*{proposition}{Proposition}
\newtheorem*{lemma}{Lemma}
\begin{document}
\maketitle
\begin{abstract} A free action of the direct product of two copies of the symmetric group on 3 elements  
on the cartesian product of two copies of the 3-sphere is
constructed.  This nonlinear action is constructed using surgery. 
The action provides a counterexample to a conjecture of Lewis made in 1968.
\end{abstract}

Results of P. A. Smith \cite{Smith (1944)} and J. Milnor \cite{Milnor (1957)} show that a dihedral group
cannot act freely on a sphere.  This implies that if a group $G$ acts freely on a sphere, then any element of
order 2 must be central.  Based on this G. Lewis \cite{Lewis (1968)} conjectured\footnote{I
thank Alejandro Adem bringing this conjecture to my attention and Peter May for providing background
for the proof.} that if a group $G$ acts freely on $S^n \times S^n$, then any subgroup isomorphic
to the Klein 4-group ($\Z_2 \times \Z_2$) must intersect the center non-trivially.  (This condition must
hold true for a free, linear action.) We will use standard techniques from surgery theory  to prove the
assertion in the title of this note and thus give a counterexample to Lewis' conjecture.  

R. G. Swan \cite{Swan (1960)} constructed a finite $CW$-complex $X$ whose fundamental group is the
symmetric group $S_3$ and whose universal cover $\widetilde{X}$ has the homotopy type of the
3-sphere $S^3$.  Milnor's result shows that $X$  does not have the homotopy type of a closed manifold. 
However, there is the following well-known result, of which we sketch a proof.

\begin{proposition}  Swan's complex $X$ is a Poincar\'e complex whose Spivak bundle \break $\nu:X \to
BG$ reduces to a stable vector bundle $\tilde{\nu} : X \to BO$.
\end{proposition}

\begin{proof} 
Swan's construction (or a spectral sequence argument) shows $H_3 \widetilde{X} \to H_3 X$  is
isomorphic to $\Z \xrightarrow{\times 6} \Z$.  Hence the transfer $H_3 X \to H_3 \widetilde{X}$ is an
isomorphism, $X$ has a fundamental class, and $X$ is a Poincar\'e complex.

A direct analysis
of the $J$-homomorphism  $\pi_*O \to \pi_*G$ shows $\pi_i(G/O) \cong \Z,0,\Z_2,0,\Z,0,\Z_2$ for $i =
0,1,2,3,4,5,6$  (see \cite{Madsen-Milgram (1979)} for background on surgery classifying
spaces).  Thus there is a single obstruction in
$H^3(X;
\pi_2(G/O))=
\Z_2$ to reduction of the Spivak bundle of $X$.  The transfer to the 3-fold cover $Y = \widetilde{X}/(12)$
 of $X$ is an injection on mod 2 cohomology, so it suffices to show that $Y$ has the homotopy type of a
closed manifold, namely $\R P^3$.  Proceed as
follows.  It is easy to find a map $f : \R P^3 \to Y$ which is an isomorphism on $\pi_1$.  Then $f$ must have
odd degree, since $f^*(a^3) = f^*(a)^3 \neq 0$, where $a \in H^1(Y;\Z_2)$ is the generator.
  By replacing $f$ by
$$h: \R P^3 \# S^3 \# \dots \# S^3  \to Y
$$
with the orientation on the $S^3$ chosen suitably and the map on the $S^3$'s given by $S^3
\tilde{\to} \widetilde{X}
\to Y$, one obtains the desired homotopy equivalence $h$.
\end{proof}

We need a computational lemma for the proof of the main theorem.

\begin{lemma}  \begin{enumerate}
\item  $\widetilde{K}_0(\Z[S_3 \times S_3]) = 0$.
\item $L_2^h(\Z[S_3 \times S_3]) = L_2^p(\F_2[S_3 \times S_3])$.
\end{enumerate}
\end{lemma}

\begin{proof}    By induction theory \cite{Dress (1975)}, it suffices to prove that $\widetilde{K}_0$
vanishes for the maximal hyperelementary subgroups: $\Z_3 \times \Z_3$ and $S_3 \times \Z_2$.  This
result for $\Z_3 \times \Z_3$ is contained in \cite[8.4]{Reiner-Ullom (1974)}.  For $S_3 \times \Z_2$ we
use the cartesian square
$$\begin{CD}
\Z[S_3 \times \Z_2] @>>> \Z[S_3] \\
@VVV @VVV \\
\Z[S_3] @>>> \F_2[S_3] 
\end{CD}$$
where the top map and left map send the generator of $\Z_2$ to $+1$ and $-1$ respectively.  Thus there
is  a Mayer-Vietoris exact sequence
\cite[Chap. 3]{Milnor (1971)}
$$K_1(\Z[S_3]) \oplus K_1(\Z[S_3]) \to K_1(\F_2[S_3]) \xrightarrow{\partial} \widetilde{K}_0(\Z[S_3
\times \Z_2]) \to \widetilde{K}_0(\Z[S_3]) \oplus \widetilde{K}_0(\Z[S_3])
$$
By \cite{Swan (1960)} or \cite{Reiner-Ullom (1974a)},  $\widetilde{K}_0(\Z[S_3])= 0$.  Now 
$\F_2[S_3] = \F_2[\Z_2] \times M_2(\F_2)$, so $K_1(\F_2[S_3])$ has order two, generated by $(12)$. 
This is in the image of the left hand map, so $\widetilde{K}_0(\Z[S_3 \times \Z_2])=0$.  Thus
$\widetilde{K}_0(\Z[S_3 \times S_3])=0$ and $L_2^h(\Z[S_3 \times S_3])=L_2^p(\Z[S_3 \times S_3])$.

The rational representation theory for $S_3 \times
S_3$ is particularly simple: $\Q[S_3 \times S_3] = \Q[S_3] \otimes_{\Q} \Q[S_3]$ is a product of
matrix rings over the rational numbers.  It follows that $L_2^p(\Z[S_3 \times S_3]) =
L_2^p(\hat{\Z}_2[S_3
\times S_3])$  where $\hat{\Z}_2$ is the 2-adic integers (see \cite[Section 4]{Hambleton-Taylor
(1998)}).  Since
$\hat{\Z}_2[S_3
\times S_3]$ is a complete semilocal ring with 2 in its Jacobson radical, $L_2^p(\hat{\Z}_2[S_3
\times S_3]) = L_2^p(\F_2[S_3
\times S_3])$ (see \cite{Wall (1973)}).

\end{proof}

\begin{theorem}  $S_3 \times S_3$ acts freely and smoothly on $S^3 \times S^3$.
\end{theorem}
\begin{proof}  By the proposition, $X$ is a Poincar\'e complex with a reducible Spivak bundle, so a
transversality argument (see \cite[Theorem 2.23]{Madsen-Milgram (1979)}) gives a degree one normal
map
$$\begin{CD}  \nu_M @>\hat{f}>> \xi \\
@VVV @VVV \\
M^3 @>f>>  X.
\end{CD}$$
There is a surgery obstruction \cite[17D]{Wall (1970)},
$$\sigma_*(f\times f,\hat{f} \times \hat{f}) \in L^h_6(\Z[S_3 \times S_3])
= L^p_6(\F_2[S_3 \times S_3])$$
which vanishes if and only if $f\times f$ is normally bordant to a homotopy equivalence.   
By the composition formula and product formula \cite{Ranicki (1980)}
\begin{align*}  
\sigma_*(f\times f,\hat{f} \times \hat{f}) &= \sigma_*(\text{Id}_M \times  f) + \sigma_*(f \times 
\text{Id}_M)\\
&= \sigma_*(f,\hat{f})\otimes (\sigma^*M + \sigma^*X) \in L^p_6(\F_2[S_3 \times S_3])
\end{align*}
with $\sigma_*(f,\hat{f}) \in
L_3^p(\F_2[S_3 \times S_3])$ and the symmetric signatures $\sigma^*M, \sigma^*X \in L_p^3(\F_2[S_3
\times S_3])$.  But $L_3^p(\F_2[S_3 \times S_3]) = 0$ since the odd projective
$L$-theory of a complete semilocal ring vanishes.  Thus the integral surgery
obstruction is zero and
$X \times X$ has the homotopy type of a  closed, smooth
manifold $N$ with fundamental group $S_3 \times S_3$ and universal cover homotopy equivalent to $S^3
\times S^3$.  

The surgery exact sequence
$$L_7(\Z) \to {\cal S}(S^3 \times S^3) \to [S^3 \times S^3,G/O] \to L_6(\Z)
$$
shows that ${\cal S}(S^3 \times S^3) = *$, since $L_7(\Z) = 0$ and $[S^3 \times S^3,G/O] \to L_6(\Z)$ is
isomorphic to $\Z_2 \xrightarrow{\text{Id}} \Z_2$, since there is a framed 6-manifold  with non-trivial 
Kervaire invariant (coincidently $S^3
\times S^3$ with the Lie invariant framing).   Thus the universal cover
of
$N$ is diffeomorphic to
$S^3 \times S^3$ and the deck transformations give the desired action.

\end{proof}

\end{document}